\documentclass[11pt,reqno]{amsart} 
\usepackage{graphicx}
\usepackage{xcolor}
\usepackage{amsmath}
\usepackage{amsfonts}
\usepackage{amssymb}
\usepackage{multicol}
\usepackage{amsthm}
\usepackage{ytableau}
\usepackage{subcaption}
\usepackage{float}
\usepackage[english]{babel}
\usepackage{hyperref}
\usepackage[capitalize,nameinlink]{cleveref}
\usepackage{stmaryrd}
\usepackage{tikz}
\usepackage[inline]{enumitem}
\colorlet{myGreen}{green!70!black}
\colorlet{darkGreen}{green!50!black}
\definecolor{myBlue}{rgb}{0.25, 0.0, 1.0}
\definecolor{lgray}{rgb}{0.75, 0.75, 0.75}
\def\ourtitle{Row Impartial Terminus}
\hypersetup{
  colorlinks=true,
  linkcolor=myBlue,
  citecolor=myGreen,
  urlcolor=myGreen,
  bookmarksopen=true,
  bookmarksnumbered,
  bookmarksopenlevel=2,
  bookmarksdepth=3
  pdftitle = {\ourtitle},
  pdfauthor= {Eric Gottlieb, Matjaz Krnc, Peter Mursic},
}
\title{\ourtitle}
\author[E. Gottlieb]{Eric Gottlieb}
\address{Rhodes College\\
2000 North Pkwy\\
Memphis, TN 38112\\
USA}
\email{gottlieb@rhodes.edu}
\author[D. Khatana]{Dawood Khatana}
\address{Rhodes College\\
2000 North Pkwy\\
Memphis, TN 38112\\
USA}
\email{dawoodkhatana15@gmail.com}
\author[M. Krnc]{Matjaž Krnc}
\address{University of Primorska\\
Titov trg 4\\
Koper, SI–6000 \\
Slovenia}
\email{matjaz.krnc@upr.si}
\author[P. Muršič]{Peter Muršič}
\address{University of Primorska\\
Titov trg 4\\
Koper, SI–6000 \\
Slovenia}
\email{peter.mursic@famnit.upr.si}
\author[I. Qureshi]{Ismael Qureshi}
\address{Rhodes College\\
2000 North Pkwy\\
Memphis, TN 38112\\
USA}
\email{mismiq64@gmail.com}
\thanks{This work was supported in part by the Slovenian Research and Innovation Agency (research projects J1-4008, J1-3002, research program P1-0383, and bilateral project BI-US-24-26-018).
The first author was supported by an internal grant from Rhodes College. 
}

\numberwithin{equation}{section}

\newcommand{\br}[1]{\llbracket #1 \rrbracket}
\newcommand{\nim}{\textsc{Nim}}
\newcommand{\rit}{\textsc{RIT}}
\newcommand{\p}{\mathcal P}
\newcommand{\n}{\mathcal N}
\newcommand{\pnim}{\textsc{PNim}}
\newcommand{\rnim}{\textsc{RNim}}
\newcommand{\rem}{\rm{rem}}
\newcommand{\core}{\rm{core}}
\newcommand{\sg}{\mathcal{G}}
\newcommand{\mex}{\rm{mex}}
\newtheorem{theorem}{Theorem}[section]
\newtheorem{lemma}[theorem]{Lemma}
\newtheorem{corollary}[theorem]{Corollary}
\newtheorem{definition}[theorem]{Definition}
\newtheorem{example}[theorem]{Example}
\newtheorem*{definition*}{Definition}
\begin{document}
\begin{abstract}
    We introduce \textsc{Row Impartial Terminus} (\rit{}), an impartial combinatorial game played on integer partitions. We show that any position in \rit{} can be uniquely decomposed into a 
    \emph{core} and a \emph{remnant}. Our central result is that the 
    Conway pair of any \rit{} position—which determines the outcome under both normal and misère play—is identical to the Conway pair of a corresponding position in the game of \nim{} defined by the remnant. 
    This finding provides a complete winning strategy for both variants of \rit{}, reducing its analysis to the well-understood framework of \nim{}. 
    As a consequence, we classify \rit{} within the Conway-Gurvich-Ho hierarchy, showing it to be forced and miserable but not pet.
\end{abstract}
\maketitle
\section{Introduction}
Several researchers have explored combinatorial games on partitions. 
In 1970, Sato \cite{sato1970maya} showed that Welter's game can be formulated this way and conjectured that its Sprague-Grundy values relate to the representation theory of the symmetric group. 
In 2018, Irie \cite{irie2018p} confirmed Sato's conjecture. 
Subsequent work by Abuku and Tada \cite{abuku2023multiple} and Motegi \cite{motegi2021gamepositionsmultiplehook} extended these results. 
Several other games on integer partitions have also been studied: \textsc{LCTR} by several authors \cite{Gottlieb2024LCTR,Ilic2019,gottlieb2025nimintegerpartitionshyperrectangles}, \textsc{CRIM} by Ba\v si\' c \cite{basic2022some,Basic_2023}, and a suite of chess-inspired impartial games, collectively called \textsc{Impartial Chess}, by Berlekamp \cite{impchess} and others \cite{gottlieb2025impartialchessintegerpartitions}. 
In her honors thesis, Meit \cite{hannah2025thesis} studied \textsc{CRPM} and \textsc{CRPS}, two partizan combinatorial games on partitions.
\subsection{Our results}
 We introduce an impartial combinatorial game on integer partitions which we call 
\textsc{Row Impartial Terminus} (\rit{}). While we focus on the standard normal play convention, our analysis extends to the misère variant as well.
We show that any \rit{} position can be uniquely decomposed into two components: a core and a remnant. The game's outcome is determined entirely by a \nim{}-position that is derived from the remnant.  
We show that from any winning position in \rit{}, there exists a winning response that is determined by playing \nim{} on the remnant. 
We prove that the 
Conway pair of any \rit{} position is identical to the \nim{} position corresponding to its remnant, providing a complete solution for both normal and misère play. 
We classify \rit{} within the Conway-Gurvich-Ho (CGH) hierarchy as forced and miserable but not pet.
\subsection{Structure of the paper}
The paper is organized as follows. 
In \cref{sec:2}, we review some definitions, notations, and results that will be used throughout the paper. 
\cref{sec:3} formally defines the rules of \rit{}, introduces the core and remnant decomposition, and establishes our main result. 
\cref{sec:cgh} extends this result to misère play and provides the game's CGH classification.
In \cref{sec:conclusion} we conclude the paper by offering directions for further work on \rit{} and related games.
\section{Preliminaries}\label{sec:2}
We denote the set of integers by $\mathbb Z$ and the set of nonnegative integers by $\mathbb Z^{\geq 0}$. For integers $i$ and $j$ we define $[i,j] = \{ n \in \mathbb Z : i \leq n \leq j\}$ and $[j]=[1,j]$. 
\subsection{Combinatorial games}
We assume the reader is familiar with basic terminology concerning combinatorial games; see Siegel \cite{Sie13}. A game consists of rules that dictate the moves the players can make from a given position.
For a game $G$, we write $p \to_{G} p'$ if there is a move from position $p$ to position $p'$, and also call such a move a $G$-move.  
Define $M_G(p)=\{p'\mid p\to_G p'\}$ to be the set of positions of $G$ reachable from $p$ in a single move. 
When clear from context, we omit the game $G$ from those notions, and also use terms \emph{position} and \emph{game} interchangeably. 
A \emph{$\p$-position} is one from which the previous player can force a win; an \emph{$\n$-position} is one from which the next player can. 
Under normal play, terminal positions are $\p$-positions and under misère play terminal positions are $\n$-positions. The other positions can be computed recursively. 
If there exist $p'$ such that $p \to p'$ and $p'$ is a $\p$-position, then $p$ is an $\n$-position. 
If $p\to p'$ implies $p'$ is an $\n$-position, then $p$ is a $\p$-position.
In impartial games under normal play, every position can be assigned a non-negative integer called a \emph{Sprague-Grundy number}. 
A position is a $\p$-position if its Sprague-Grundy number is zero, and is an $\n$-position otherwise.
Let $S$ be a proper subset of the nonnegative integers. 
The \emph{mex}, or minimal excludant, of $S$ is defined to be the smallest nonnegative integer not in $S$. For example,
\[
\mex(\{0,1, 2, 4\}) = 3 \quad \mbox{and} \quad \mex(\{1, 2, 3\}) = 0.
\]
The \emph{Sprague-Grundy value} of a position $x$ in game is recursively defined by
\[
\sg_G(x) = \mex\bigl(\{\sg(y) \mid x\to y\}\bigr).
\]
Note that $\sg(p) = 0$ for every terminal position $p$. The \emph{misère Grundy} value $\sg_G^-(p)$ is defined the same way as $\sg_G(p)$, except terminals are assigned value $1$. 
The \emph{Conway pair} of a position $p$ is defined to be $(\sg_G(p),\sg_G^-(p))$, and we say that $p$ is an $(\sg_G(p),\sg_G^-(p))$-position.
Sprague-Grundy values are useful for analyzing the \emph{disjunctive sum} $G_1 + \dots + G_k$ of normal-form impartial games games $G_1, \dots, G_k$. In $G_1 + \dots + G_k$, a position consists of a $k$-tuple $(p_1, \dots, p_k)$, where $p_j$ is a position of $G_j$. Moves in $G_1 + \dots + G_k$ are of the form $(p_1, \dots, p_k) \to (\bar p_1, \dots, \bar p_k)$, where $p_j \to \bar p_j$ in $G_j$ and $p_i= \bar p_i$ for $i \neq j$. A terminal position is one of the form $(p_1, \dots, p_k)$ where each $p_j$ is terminal in $G_j$. 
For integers $a$ and $b$, their \emph{nim-sum} $a\oplus b$ is the bitwise XOR of their binary representations. 
For example,
\[
5 \oplus 3 = 101_2 \oplus 011_2 = 110_2 = 6.
\]
Nim-sum is associative, commutative, and satisfies $a \oplus a = 0$. 
Furthermore, a theorem of Sprague \cite{Spr35} and Grundy \cite{Gru39} asserts that 
\[\sg_{G_1 + \dots + G_k}((p_1, \dots, p_k)) = \sg_{G_1}(p_1) \oplus \dots \oplus \sg_{G_k}(p_k).\] 
Thus, $(p_1, \dots, p_k)$ is a $\p$-position in $G_1 + \dots + G_k$ if and only if $\sg(p_1) \oplus \dots \oplus \sg(p_k) = 0.$
\nim{} is a classic combinatorial game foundational to the Sprague-Grundy theorem. It is played with heaps of objects; on their turn, a player removes any number of objects from a single heap.
\begin{definition}\label{ex:NIM}
The game $1$-\nim{} is played on  $\mathbb Z_{\ge 0}$ and for $p,p'\in\mathbb Z_{\ge 0}$, $p\to p'$ is a move of $1$-\nim{} whenever $p-p'>0$.
$\nim{}$ is defined to be the disjoint sum $G_1+\dots+G_k$ for any positive integer $k$, where each $G_i$ is an instance of $1$-\nim{}.
\begin{example} 
The $\nim{}$-position $(3, 1, 2)$ allows six moves; one of them is to $(1,1,2)$. 
Furthermore, $(3, 1, 2)$ is a $\p$-position since $3 \oplus 1 \oplus 2 = 0$. 
On the other hand, $(1,1,2)$ is an $\n$-position since $1 \oplus 1 \oplus 2 = 2$. 
\end{example}
\end{definition}
\subsection{Conway-Gurvich-Ho classification} \label{sec:cgh}
The interaction of the normal and mis\`ere forms of combinatorial games was initiated by Conway \cite{Con76} and expanded by Gurvich and Ho \cite{GURVICH201854}. The so-called \emph{CGH classification} is defined by certain properties which we now recall.
The domestic class is not relevant for this paper so we omit it.
\begin{definition*}[CGH classification] \label{D.DTP}
An impartial game is called \begin{enumerate} 
\item {\em forced} if each move from a $(0,1)$-position results in a $(1,0)$-position and vice versa;
\item  {\em miserable} if for every position $x$, one of the following holds: (i) $x$ is a $(0, 1)$-position or $(1, 0)$-position, or (ii) there is no move from $x$ to a $(0,1)$-position or $(1, 0)$-position, or (iii) there are moves from $x$ to both a $(0, 1)$-position and a $(1, 0)$-position;
\item {\em pet} if it has only $(0,1)$-positions, $(1,0)$-positions, and $(k,k)$-positions with $k \geq 2$;
\end{enumerate}
\end{definition*}
According to the above definitions, tame, pet, miserable, returnable and forced games form nested classes:
a pet game is returnable and miserable, a miserable game is tame, 
and a forced game is returnable.
\subsection{Partitions}
Let $n$ be a nonnegative integer. An \emph{(integer) partition of $n$} is a way to write $n$ as a sum of positive integers in nonincreasing order. If $\lambda_1 + \dots + \lambda_r$ is a partition of $n$, then we write $\lambda = \br{\lambda_1, \dots, \lambda_r}$ and $\lambda \vdash n$ and refer to the $\lambda_j$'s as the \emph{parts} of $\lambda$. The only partition of zero is the empty partition $\br{\, }$. 
The partition $\lambda = \br{\lambda_1, \dots, \lambda_r}$ can be represented using a \emph{Young diagram}. 
This is a left-justified array of boxes in which there are $\lambda_j$ boxes in the $j$th row from the top.
In this paper, we do not distinguish between $\lambda$ and its Young diagram. 
The following Figure 1 shows the Young diagrams for the partitions of 4:
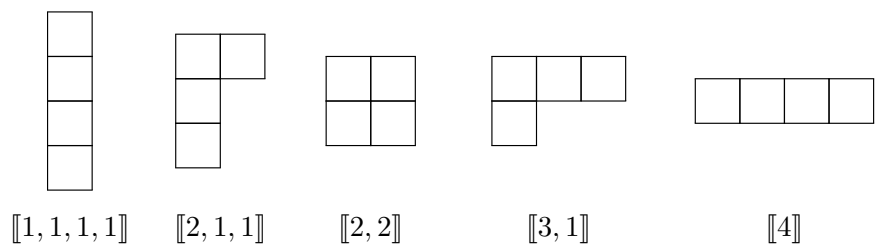
\begin{figure}
\begin{tikzpicture}
\node[scale=1] (a) at (0,0){\begin{ytableau}
    ~ \\
    ~ \\
    ~ \\
    ~
\end{ytableau}};
\node[scale=1] (b) at (2,0){\begin{ytableau}
    ~ & ~ \\
    ~ \\
    ~
\end{ytableau}};
\node[scale=1] (c) at (4,0){\begin{ytableau}
    ~ & ~ \\
    ~ & ~
\end{ytableau}};
\node[scale=1] (d) at (6.5,0){\begin{ytableau}
    ~ & ~ & ~ \\
    ~ 
\end{ytableau}};
\node[scale=1] (e) at (9.5,0){\begin{ytableau}
    ~ & ~ & ~ & ~
\end{ytableau}};
\node at (0,-1.7){$\br{1,1,1,1}$};
\node at (2,-1.7){$\br{2,1,1}$};
\node at (4,-1.7){$\br{2,2}$};
\node at (6.5,-1.7){$\br{3,1}$};
\node at (9.5,-1.7){$\br{4}$};
\end{tikzpicture}
\caption{The Young diagrams of the  partitions of $4$.}
\end{figure}
\section{Row Impartial Terminus}\label{sec:3}
\begin{figure}[H]
    \centering
    \ytableausetup{boxsize=2em}
    \begin{tikzpicture}[scale=1.2]
    \node[scale=0.6] (a) at (-1,0){    \begin{ytableau}
    ~ & ~ & ~ & ~ & 5 \\ 
    ~ & ~ & 3 & 4 \\ 
    ~ & 2 \\   
    1
    \end{ytableau}};
    \node[scale=0.6] (a1) at (5,0){    \begin{ytableau}
        ~ & ~ & ~ & ~ \\ 
        ~ & ~ & ~ & ~ \\ 
        ~ & ~ \\   
        ~
    \end{ytableau}};
        \node[scale=0.6] (a2) at (5,-2.5){    \begin{ytableau}
        ~ & ~ & ~ & ~ & ~ \\ 
        ~ & ~ & ~ \\ 
        ~ & ~ \\   
        ~
    \end{ytableau}};
        \node[scale=0.6] (a3) at (5,-5){    \begin{ytableau}
        ~ & ~ & ~ & ~ & ~ \\ 
        ~ & ~  \\ 
        ~ & ~ \\   
        ~
    \end{ytableau}};
        \node[scale=0.6] (a4) at (2,-5){    \begin{ytableau}
        ~ & ~ & ~ & ~ & ~ \\ 
        ~ & ~  & ~ & ~ \\ 
        ~  \\   
        ~
    \end{ytableau}};
        \node[scale=0.6] (a5) at (-1,-5){    \begin{ytableau}
        ~ & ~ & ~ & ~ & ~ \\ 
        ~ & ~ & ~ & ~ \\ 
        ~ & ~ \\  
    \end{ytableau}};
    \draw[shorten >=0.5cm,shorten <=1cm,->,>=stealth, thick] (a)--(a1)node[midway,above,scale=0.7]{5};
    \draw[shorten >=0.5cm,shorten <=1cm,->,>=stealth, thick]  (a)--(a2)node[midway,above,shift={(0.3,-0.1)},scale=0.7]{4};
    \draw[shorten >=0.5cm,shorten <=1cm,->, >=stealth, thick] (a)--(a3)node[midway,above right,scale=0.7]{3};
    \draw[shorten >=0.5cm,shorten <=1cm,->,>=stealth, thick]  (a)--(a4)node[midway,right,scale=0.7]{2};
    \draw[shorten >=0.5cm,shorten <=1cm,->,>=stealth, thick]  (a)--(a5)node[midway,right,scale=0.7]{1};
    \end{tikzpicture}
    \caption{\label{fig:ydmoves}Young diagram of $\br{5,4,2,1}$ along with all possible moves.
    }
    
\end{figure}
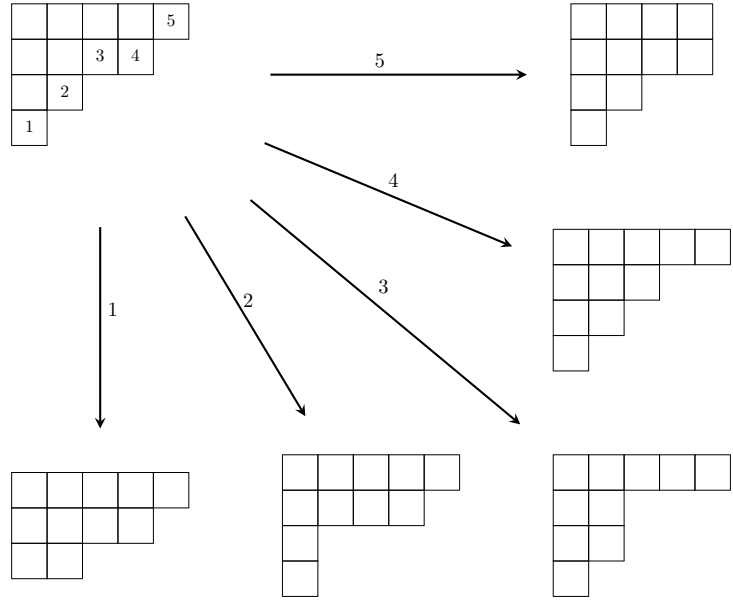

\textsc{Row Impartial Terminus}, or \rit{}, is an impartial combinatorial game played on Young diagrams. 
For any partition $\lambda=\br{\lambda_1, \dots, \lambda_r}$ and $k\in [\lambda_1]$, 
there is a move from $\lambda$ to  $\bar\lambda=\br{\bar\lambda_1, \dots, \bar\lambda_r}$ where 
$\bar\lambda_i=k-1$
if $i$ is the largest integer such that $\lambda_i\ge k$, and
$\bar\lambda_j=\lambda_j$ for $j\neq i$.
We will refer to the move above as a \emph{move on $i$th row} of $\lambda$.
The only terminal position is the empty partition. 
Less formally, a move consists of shortening one of the rows of $\lambda$ so that the result remains a Young diagram; see \cref{fig:ydmoves}.
Given a partition $\lambda = \br{\lambda_1, \ldots, \lambda_r}$ define $\core(\lambda)$ to be the partition \[\br{\lambda_2, \lambda_2, \lambda_4, \lambda_4, \ldots, \lambda_{2\lfloor r/2 \rfloor}, \lambda_{2\lfloor r/2 \rfloor}}\] and $\rem(\lambda)$ to be the $\lceil r/2 \rceil$-tuple 
\[(\lambda_1 - \lambda_2,
\lambda_3 - \lambda_4,
\ldots,  \lambda_{2\lceil r/2 \rceil-1}-\lambda_{2\lceil r/2 \rceil}),\] where $\lambda_j = 0$ if $j > r$. 
\begin{lemma}\label{lem:oddtonim}
    Let $\lambda$ be any non-empty partition and suppose that $\bar \lambda$ is the result of a move on an odd-numbered row of $\lambda$. 
    Then $\core(\bar\lambda) = \core(\lambda)$ and $\rem(\lambda)\to_{\nim}\rem(\bar \lambda)$. 
    
    Conversely, for any \nim{}-move  $\rem(\lambda)\to p$ there exists a \rit{}-move $\lambda \to\bar \lambda$ on an odd-numbered row for which $\core(\bar \lambda) = \core(\lambda)$ and $\rem(\bar\lambda)=p$. 
\end{lemma}
\begin{proof}
Let $\lambda = \br{\lambda_1, \ldots, \lambda_r}$ and $\bar\lambda = \br{\bar\lambda_1, \ldots, \bar\lambda_r}$. By the first part of the statement we have that $\lambda_{j}-\bar\lambda_{j}=i$ for some positive $i$ and odd $j\in [r]$ and $\lambda_k=\bar\lambda_k$ for $k\neq j$. 
By definition we have
\begin{align*}
\rem(\lambda) &= \left(\lambda_1 - \lambda_2, \ \ldots\ , 
\lambda_{j} - \lambda_{j+1}, \ \phantom{-i}  \ \ldots\ ,\lambda_{2\lceil r/2 \rceil-1}-\lambda_{2\lceil r/2 \rceil} \right)\\
\rem(\bar\lambda) &= \left(\lambda_1 - \lambda_2, \ \ldots\ , 
\lambda_{j} - \lambda_{j+1}-i,  \ \ldots\ ,\lambda_{2\lceil r/2 \rceil-1}- \lambda_{2\lceil r/2 \rceil} \right),
\end{align*}
in particular, $\rem(\lambda)\to \rem(\bar\lambda)$ is a \nim{}-move. To conclude the first part of the claim, observe that since \(\text{core}(\lambda)\) depends only on the even-indexed elements of \(\lambda\), a move on an odd-numbered row does not affect  \(\text{core}(\lambda)\), thus $\core(\lambda)=\core(\bar\lambda)$.
Conversely, a \nim{}-move in the $i$-th coordinate of $\rem(\lambda)$ produces 
\begin{align}
\tau=\left(\lambda_1 - \lambda_2, \ldots, \lambda_{2i-1} - \lambda_{2i}-k, \ldots,\lambda_{2\lceil r/2 \rceil-1}, \lambda_{2\lceil r/2 \rceil} \right),\label{eq:tau}
\end{align}
where $1 \leq k \leq \lambda_{2i-1} - \lambda_{2i}$.
Since $\lambda_{2i-1} - \lambda_{2i}-k\ge 0$, there is a \rit{}-move on the $(2i-1)$-th row of $\lambda$ leading to the  position 
\begin{align}
\label{eq:rit-position}
\bar\lambda = \br{\lambda_1, \ldots, \lambda_{2i-1} - k, \lambda_{2i}, \ldots, \lambda_r}.
\end{align}
By definition of $\rem$,  (\ref{eq:tau}) and (\ref{eq:rit-position}) imply $\rem{(\bar\lambda)}=\tau$.
\end{proof}
\begin{lemma}\label{lem:two-move-rem}
Let \(\lambda = \br{\lambda_1, \dots, \lambda_r}\) be a partition with $r \geq 2$. Let \(\overline{\lambda}\) be the result of a \rit{}-move on an even row of \(\lambda\). Then, there exists a \rit{}-move on an odd row of \(\overline{\lambda}\), resulting in a partition \(\tilde{\lambda}\), such that \(\rem(\tilde{\lambda}) = \rem(\lambda)\).
\end{lemma}
\begin{proof}
Let $\bar\lambda = \br{\bar\lambda_1, \dots, \bar\lambda_r}$, and let $j$ be an even integer such that $\lambda\to_{\rit{}}\bar\lambda$ and $\lambda_{j}-\bar\lambda_{j}=i$ for some positive $i$. 
We claim that the partition 
\[\bar {\bar \lambda}=\br{\bar\lambda_1, \dots,\bar\lambda_{j-1},\bar\lambda_{j}-i,\bar\lambda_{j+1}, \dots,\bar\lambda_r}\]
satisfies $\rem(\bar{\bar \lambda}) = \rem(\lambda)$. To see this, note that we have $\lambda_{j-1} - \lambda_{j} \geq 0$, so $\bar \lambda_{j-1} - \bar \lambda_{j} = \lambda_{j-1} - (\lambda_{j} -i) \geq i$. Thus
\begin{align*}
    \rem(\bar{\bar \lambda}) &= \left(\bar{\bar \lambda}_1 - \bar{\bar \lambda}_2, \ldots, \bar{\bar \lambda}_{j-1}-\bar{\bar \lambda}_{j}, \ldots, \bar{\bar \lambda}_{2 \lceil r/2 \rceil -1} - \bar{\bar \lambda}_{2 \lceil r/2 \rceil} \right) \\
    &= \left(\lambda_1 - \lambda_2, \ldots, (\lambda_{j-1} - i) - (\lambda_{j} -i), \ldots, \lambda_{2 \lceil r/2 \rceil -1} - \lambda_{2 \lceil r/2 \rceil} \right) \\
    &= \left(\lambda_1 - \lambda_2, \ldots, \lambda_{j-1} - \lambda_{j}, \ldots, \lambda_{2 \lceil r/2 \rceil -1} - \lambda_{2 \lceil r/2 \rceil} \right) \\
    &= \rem(\lambda). \qedhere
\end{align*}
\end{proof}
\begin{theorem}\label{thm:grit=grem}
    For any partition $\lambda$ we have $\sg_{\rit}(\lambda) = \sg_{\nim}(\rem(\lambda))$.
\end{theorem}
\begin{proof}
Let \(\lambda = \br{\lambda_1, \dots, \lambda_r}\) be a partition of $n$.
We proceed by induction on $n$.
Since $\sg_{\rit}(\br{\,}) = 0$ and  $\rem(\br{\,}) = (\,)$ and  $\sg_{\nim}((\,)) = 0$, this settles the base case $n=0$.
    
For the induction step,   
    assume that for  all partitions in ${\bar\lambda}$ such that $\lambda\to \bar\lambda$ we have $\sg_{\rit}(\bar \lambda) = \sg_{\nim}(\rem(\bar \lambda))$. We want to show that
    \[
    \sg_{\rit}(\lambda) = \sg_{\nim}(\rem(\lambda)).
    \]
This is equivalent to showing that
\[
\underset{\bar{\lambda}\in M_\rit(\lambda)}{\mex}\left(\sg_{\rit}({\bar\lambda})\right) = 
\underset{\tau\in M_\nim(\rem(\lambda))}{\mex}
\left(\sg_{\nim}(\tau)\right).
\]
    \(M_{\rit}(\lambda)\) can be divided into 
    \(M_{\rit}(\lambda) = M_{\rit}^{\text{odd}}(\lambda) \cup M_{\rit}^{\text{even}}(\lambda)\), 
    where the odd and even superscripts denote a move in either an odd or even row of $\lambda$. 
     Due to \cref{lem:oddtonim}, the restricted function $\rem:M_{\rit}^{\text{odd}}(\lambda) \to M_{\nim}(\rem(\lambda))$ is a bijection. 
    By the induction hypothesis, for any $\overline{\lambda}\in M_\rit(\lambda)$ we have $\sg(\overline{\lambda}) = \sg(\rem(\overline{\lambda}))$, so
    \[
     \{\sg_\nim({\tilde{\lambda}}) \mid \tilde{\lambda} \in M_{\nim}(\rem(\lambda))\}=\{\sg_\rit({\overline{\lambda}}) \mid \overline{\lambda} \in M_{\rit}^{\text{odd}}(\lambda)\},
    \]
    which implies
    \[
    \sg_\nim(\rem(\lambda)) = \mex(\{\sg_\rit({\overline{\lambda}}) \mid \overline{\lambda} \in M_{\rit}^{\text{odd}}(\lambda)\}).
    \]
    To conclude the proof it is enough to show that, for any $\overline{\lambda} \in M_{\rit}^{\text{even}}(\lambda)$, we have $\sg_\rit(\overline{\lambda}) \neq\sg_\nim(\rem(\lambda))$.
    To this end, let  $\overline{\lambda} \in M_{\rit}^{\text{even}}(\lambda)$ where $\overline{\lambda} = \left(\lambda_1,\dots,\lambda_{2i} - k,\dots,\lambda_r\right)$ is the result of a move removing $k$ blocks from row $2i$ of $\lambda$. 
    As a response, removing $k$ blocks from the $(2i-1)$th row of $\overline{\lambda}$ produces $\tilde{\lambda} = \left(\lambda_1,\dots,\lambda_{2i-1}-k,\lambda_{2i} - k,\dots,\lambda_r\right)$. 
    Due to \cref{lem:two-move-rem}, we have $\rem(\tilde{\lambda})=\rem(\lambda)$. Finally,  $\overline{\lambda}\to_{\rit}\tilde{\lambda}$ implies that $\sg(\overline{\lambda}) \neq\sg(\rem(\lambda))$, as desired. 
\end{proof}
\begin{corollary}\label{cor:howtoplay}
Let $\lambda$ be a losing position under \rit{}. Regardless of the opponents move, there exists a winning response by reducing an odd-indexed row.
Furthermore, if the opponent reduces an even indexed row $i$, then the winning response consists of reducing $(i-1)$th row by the same amount.
\end{corollary}
By \cref{thm:grit=grem,cor:howtoplay}, we know that the $\sg$ of the \rit{} game is equal to the $\sg$ of $\rem(\lambda)$ played in \nim{}. The objective should be to always reduce the $\sg$ value of $\rem(\lambda)$ to zero, by choosing the correct odd-indexed row. 
\section{Conway-Gurvich-Ho classification}\label{sec:cgh}
\cref{thm:grit=grem} shows that the Sprague-Grundy value of any position in \rit{} corresponds to the Sprague-Grundy value of its remnant (under \nim). In this section we show that this results extends to Conway pairs, which allows us to place \rit{} in the Conway-Gurvich-Ho classification scheme.
\begin{theorem}\label{thm:misere-grit=grem}
    For any partition $\lambda$ we have $\sg^-_{\rit}(\lambda) = \sg^-_{\nim}(\rem(\lambda))$.
\end{theorem}
\begin{proof}
Let \(\lambda = \br{\lambda_1, \dots, \lambda_r}\) be a partition of $n$.
We proceed by induction on $n$.
Since $\sg^-_{\rit}(\br{\,}) = 1$ and  $\rem(\br{\,}) = (\,)$ and  $\sg^-_{\nim}((\,)) = 1$, this settles the base case $n=0$.
We omit the induction step, as it is exactly the same as the one from the proof of \cref{thm:grit=grem}, modulo replacing $\sg$ with $\sg^-$.
\end{proof}
\begin{corollary}\label{cor:cpair}
    The Conway pair of any partition $\lambda$ under \rit{} is equal to the Conway pair of $\rem{(\lambda)}$ under \nim. 
    In particular, \rit{} is forced and miserable, but not pet.
\end{corollary}
\begin{figure}
    \centering
    \begin{tikzpicture}[scale=0.5]
\def\unknown{unknown}
\def\textsizes{0.6}
\draw [darkGreen] (2,4.5)--node[fill=none,yshift=1mm,xshift=12mm,scale=\textsizes]{miserable}++(22,0)--++(0,7.5)--++(-22,0)--cycle;
\draw [red] (3,7.5)--node[fill=none,yshift=1mm,xshift=22mm,scale=\textsizes]{pet}++(16,0)--++(0,3.5)--++(-16,0)--cycle;
\draw [gray] (3.5,1)--node[fill=none,yshift=1mm,scale=\textsizes]{forced}++(8,0)--++(0,9.5)--++(-8,0)--cycle;
\draw [cyan] (2.5,0.5)--node[fill=none,yshift=1mm,xshift=-17.5mm,,scale=\textsizes]{returnable}++(17,0)--++(0,11)--++(-17,0)--cycle;
\draw [blue] (1.5,1.5)--node[fill=none,yshift=1mm,xshift=12mm,scale=\textsizes]{tame}++(23,0)--++(0,11)--++(-23,0)--cycle;
\node at (15.27,9.5) {$1$-\pnim{}, $1$-\rnim{},};
\node at (15.27,8.5) {\textsc{Subtraction}};
\node at (7.5,9) {$1$-\nim{}, $2$-RIT};
\node at (7.5,6) {\nim{}, \rit{}};
\node at (15.3,6.5) { \pnim{}, \rnim{},};
\node at (15.3,5.5) {$\textsc{Wythoff}$};
\node at (15.45,3.5){\textsc{Downright},};
\node at (15.45,2.5){\textsc{LCTR}};
\node at (22,1) {\textsc{Rook}};
\end{tikzpicture}
    \caption{The Conway-Gurvich-Ho classification of \rit{} alongside other impartial games, where $k$-\rit{} is \rit{} restricted to partitions with at most $k$ rows.}
    \label{fig:conwaygurvich}
\end{figure}
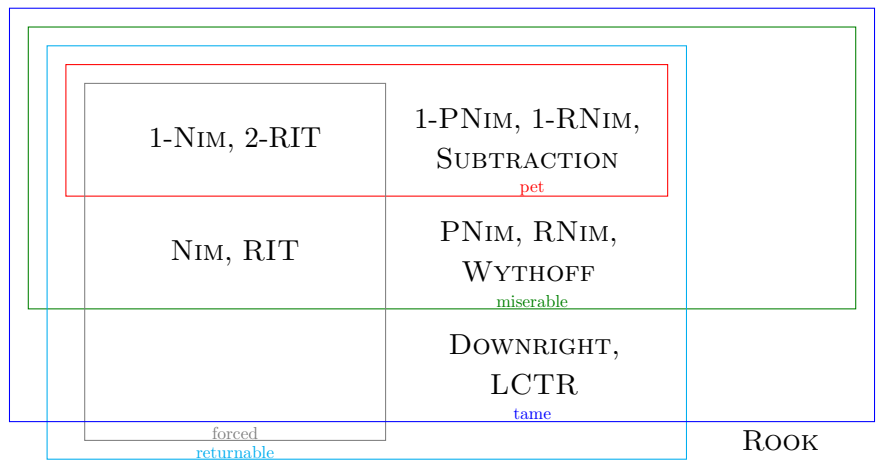
\section{Conclusion}\label{sec:conclusion}
This paper introduces \rit{}, a new impartial game on partitions in which players take turns removing parts of rows while maintaining a valid Young diagram. We divide the Young diagram into two parts: the core and the remnant. 
We establish that the $\sg$ of the remnant played in \nim{} is equal to the game of \rit{}. We analyze different possible moves and state the order of moves necessary to win this game. 
\cref{thm:grit=grem} directly implies the $\p$ or $\n$ status for several known families, such as staircases or quadrated partitions. It would be interesting to systemmatically classify other partition families, and also to explore more derivations of the game branched from \rit{}.
\bibliography{bibligraphy}
\bibliographystyle{plain}
\end{document}